\documentclass[11pt,dvipdfm]{article}
\usepackage{amsfonts}
\usepackage{graphicx}
\usepackage{amsfonts, amsmath, amssymb, latexsym, epsfig}
\usepackage{mathrsfs}
\newcommand{\blst}{\begin{trivlist}}
\newcommand{\elst}{\end{trivlist}}

\newtheorem{theorem}{\hspace{1.3em}Theorem}[section]

\newtheorem{lemma}{\hspace{1.3em}Lemma}[section]
\newtheorem{definition}{\hspace{1.3em}Definition}[section]

\newtheorem{thm}{Theorem}[section]
\newtheorem{prop}[thm]{Proposition}
\newtheorem{cor}[thm]{Corollary}
\newtheorem{lem}[thm]{Lemma}
\newtheorem{conj}[thm]{Conjecture}
\newtheorem{exa}[thm]{Example}
\newtheorem{defn}[thm]{Definition}

\newcommand{\ben}{\begin{enumerate}}
\newcommand{\een}{\end{enumerate}}
\newcommand{\ble}{\begin{lem}}
\newcommand{\ele}{\end{lem}}
\newcommand{\bth}{\begin{thm}}
\renewcommand{\eth}{\end{thm}}
\newcommand{\bpr}{\begin{prop}}
\newcommand{\epr}{\end{prop}}
\newcommand{\bco}{\begin{cor}}
\newcommand{\eco}{\end{cor}}
\newcommand{\bcon}{\begin{conj}}
\newcommand{\econ}{\end{conj}}
\newcommand{\bde}{\begin{defn}}
\newcommand{\ede}{\end{defn}}
\newcommand{\bex}{\begin{exa}}
\newcommand{\eex}{\end{exa}}
\newcommand{\barr}{\begin{array}}
\newcommand{\earr}{\end{array}}
\newcommand{\btab}{\begin{tabular}}
\newcommand{\etab}{\end{tabular}}
\newcommand{\beq}{\begin{equation}}
\newcommand{\eeq}{\end{equation}}

\newcommand{\bea}{\begin{eqnarray*}}
\newcommand{\eea}{\end{eqnarray*}}
\newcommand{\beaa}{\begin{eqnarray}}
\newcommand{\eeaa}{\end{eqnarray}}
\newcommand{\bce}{\begin{center}}
\newcommand{\ece}{\end{center}}
\newcommand{\bpi}{\begin{picture}}
\newcommand{\epi}{\end{picture}}
\newcommand{\bfi}{\begin{figure} \begin{center}}
\newcommand{\efi}{\end{center} \end{figure}}

\newcommand{\bsl}{\begin{slide}{}}
\newcommand{\esl}{\end{slide}}

{}
{}
{}
{}
{}
{}
{}
{}
{}
{}
{}
{}
{}
{}
{}
{}
{}
{}
{}
{}
{}
{}
{}
{}
{}
{}
{}
{}
\newenvironment{proof}{
\par
\noindent {\bf Proof.}\rm}{\mbox{}\hfill\rule{0.5em}{0.809em}\par}
\begin{document}

\title{$(G,m)$-Multiparking Functions}
\author{Hungyung Chang$^{a}$
\and Po-Yi Huang$^{b,}$\thanks{Partially supported by NSC
96-2115-M-006-012 }
 \and Jun Ma$^{c,}$\thanks{Email address of the corresponding author: majun@math.sinica.edu.tw}\\
 \and Yeong-Nan Yeh$^{d,}$\thanks{Partially supported by NSC 96-2115-M-001-005}
}\date{} \maketitle \vspace*{-1.2cm}\begin{center} \footnotesize
$^{a,c,d}$ Institute of Mathematics, Academia Sinica, Taipei, Taiwan\\
$^{b}$ Department of Mathematics, National Cheng Kung University, Tainan, Taiwan\\

\end{center}

 \vspace*{-0.3cm}
\thispagestyle{empty}
\begin{abstract}
The conceptions of $G$-parking functions and $G$-multiparking
functions were introduced in \cite{postnikov2004} and \cite{kostic}
respectively. In this paper, let $G$ be a connected graph with
vertex set $\{1,2,\ldots,n\}$ and $m\in V(G)$. We give the
definition of $(G,m)$-multiparking function. This definition unifies
the conceptions of $G$-parking function and $G$-multiparking
function. We construct bijections between the set of
$(G,m)$-multiparking functions and the set of $\mathcal{F}_{G,m}$ of
spanning color $m$-forests of $G$. Furthermore we define  the
$(G,m)$-multiparking complement function, give the reciprocity
theorem for $(G,m)$-multiparking function and extend the results
\cite{Y3,kostic} to $(G,m)$-multiparking function. Finally, we use a
combinatorial methods to give a recursion of the generating function
of the sum $\sum\limits_{i=1}^na_i$ of $G$-parking functions
$(a_1,\ldots,a_n)$.
\end{abstract}

\section{Introduction}
J. Riordan \cite{R} define the parking function as follows: $m$
parking spaces are arranged in a line, numbered $1$ to $n$ left to
right; $n$ cars, arriving successively, have initial parking
preferences, $a_i$ for $i$, chosen independently and at random;
$(a_1,\cdots,a_n)$ is called preference function; if space $a_i$ is
occupied, car $i$ moves to the first unoccupied space to the right;
if all the cars can be parked, then  the preference function is
called parking function.

Konheim and Weiss \cite{konhein1966} introduced the conception of
the  parking functions of length $n$ in the study of the linear
probes of random hashing function.  J. Riordan \cite{R} studied  the
parking functions  and derived that the number of parking functions
of length $n$ is $(n+1)^{n-1}$, which coincides with the number of
labeled trees on $n+1$ vertices by Cayley's formula. Several
bijections between the two sets are known (e.g., see
\cite{FR,R,SMP}). Parking functions have been found in connection to
many other combinatorial structures such as acyclic mappings,
polytopes, non-crossing partitions, non-nesting partitions,
hyperplane arrangements, etc. Refer to \cite{F,FR,GK,PS,SRP,SRP2}
for more information.

Parking function $(a_1,\cdots,a_n)$ can be redefined that its
increasing rearrangement $(b_1,\cdots,b_n)$ satisfies $b_i\leq i$.
 Pitman and  Stanley generalized the notion of parking functions
in \cite{PS}. Let ${\bf x}=(x_1,\cdots,x_n)$ be a sequence of
positive integers. The sequence $\alpha=(a_1,\cdots,a_n)$  is called
an ${\bf x}$-parking function if the non-decreasing rearrangement
$(b_1,\cdots,b_n)$ of $\alpha$ satisfies $b_i\leq x_1+\cdots +x_i$
for any $1\leq i\leq n$. Thus, the ordinary parking function is the
case ${\bf x}=(1,\cdots,1)$. By the determinant formula of
Gon\v{c}arove polynomials, Kung and Yan \cite{KY} obtained the
number of ${\bf x}$-parking functions for an arbitrary ${\bf x}$.
See also \cite{Y1,Y2,Y3} for the explicit formulas and properties
for some specified cases of ${\bf x}$.

Recently, Postnikov and Shapiro \cite{postnikov2004} gave a new
generalization, building on work of Cori, Rossin and Salvy
\cite{cori2002}, the $G$-parking functions of a graph. For the
complete graph $G=K_{n+1}$, the defined functions in
\cite{postnikov2004} are exactly the classical parking functions.
Chebikin and Pylyavskyy \cite{Denis2005} established a family of
bijections from the set of $G$-parking functions to the spanning
trees of that graph. Then Dimitrije Kostic and Catherine H. Yan
\cite{kostic} proposed the notion of a $G$-multiparking function, a
natural extension of the notion of a $G$-parking function and
extended the result of \cite{Y3} to arbitrary graphs. They
constructed a family of bijections from the set of $G$-multiparking
functions to the spanning forests of $G$. By the definition in
\cite{kostic}, it is easy to see that the vertex $1$ is always the
root in all $G$-multiparking functions $f$ if the vertex set of $G$
is $\{1,2,\ldots,n\}$. One of the motivations of this paper is to
consider the case in which the vertex $1$ isn't the root. So, we
give the definition of $(G,m)$-multiparking function. This
definition unifies the conceptions of $G$-parking function and
$G$-multiparking function because a $(G,m)$-multiparking function is
a $G$-multiparking functions and a $G$-parking functions when $m=1$
and $m=n$ respectively. Using the methods developed by Dimitrije
Kostic and Catherine H. Yan \cite{kostic}, we construct bijections
between the set $\mathcal{MP}_{G,m}$ of $(G,m)$-multiparking
functions and the set of $\mathcal{F}_{G,m}$ of spanning color
$m$-forests of $G$.

Richard Stanley's book \cite{SRP4}, in the context of rational
generating functions, devotes an entire section to exploring the
relationships (called reciprocity relationships) between positively-
and nonpositively-indexed terms of a sequence. If
$\alpha=(a_1,a_2,\cdots,a_n)$ is a $K_{n+1}$-parking function of
length $n$, then $(n-a_1,n-a_2,\cdots,n-a_n)$ is called a complement
of the parking function $\alpha$. It is easy to see that the sums
$\sum\limits_{i=1}^{n}(n-a_i)-{n+1\choose{2}}$ and
$\sum\limits_{i=1}^na_i$ are connected with the reciprocity law for
$K_{n+1}$-parking function. Also the sum
$\sum\limits_{i=1}^{n}(n-a_i)-{n+1\choose{2}}$ is one of the most
important statistic of $K_{n+1}$-parking function. Knuth
\cite{Knuth1998} indicated that it corresponds to the number of
linear probes in hashing functions. Kreweras \cite{Kreweras}
concluded that it is equal to the number of inversions in labeled
trees on $[n+1]$. Also it coincides with the number of hyperplanes
separating a given region from the base region in the extended Shi
arrangements \cite{SRP3}. Catherine H. Yan \cite{Y3} gave a
combinatorial explanation, which revealed the underlying
correspondence between the classical parking functions and labeled,
connected graphs. Furthermore, for arbitrary graph $G$, Dimitrije
Kostic and Catherine H. Yan \cite{kostic} indicate the relations
between $G$-inversions and $G$-multiparking functions. There is a
interesting problem: how to define the complement of a
$(G,m)$-multiparking function? In this paper, we define  the
$(G,m)$-multiparking complement function, give the reciprocity
theorem for $(G,m)$-multiparking function and extend the results
\cite{Y3,kostic} to $(G,m)$-multiparking function.

In \cite{kostic}, Dimitrije Kostic and Catherine H. Yan related
$G$-multiparking functions to the Tutte polynomial $T_G(x,y)$ of
$G$. The generating function $P_{G}(q)$ of $G$-parking functions is
defined as $P_{G}(q)=\sum\limits_{f}q^{\sum\limits_{i\in V(G)}f(i)}$
where $f$ ranges over all $G$-parking functions. Their results
implies that $P_{G}(q)=q^{|E(G)|-|V(G)|}T_G(1,\frac{1}{q})$. Note
that the Tutte polynomial of $G$ satisfies the recursion
$$T_G(x,y)=\left\{\begin{array}{lll}
xT_{G-e}(x,y)&\text{if}& e\text{ is a bridge}\\
yT_{G-e}(x,y)&\text{if}& e\text{ is a loop}\\
T_{G-e}(x,y)+T_{G\setminus e}(x,y)&\text{otherwise}&
\end{array}\right.$$ where $G-e$ is a graph
obtained by deleting the edge $e$ and $G{\setminus e}$ is a graph
obtained from $G$ contracting the the vertices $i$ and $j$ from $G$.
By a combinatorial method, we give a recursion of the generating
function of the sum $\sum\limits_{i=1}^na_i$ of $G$-parking
functions $(a_1,\ldots,a_n)$.

This paper is organized as follows. In Section $2$, we construct
bijections between the set $\mathcal{MP}_{G,m}$ of
$(G,m)$-multiparking functions and the set of $\mathcal{F}_{G,m}$ of
spanning color $m$-forests of $G$. In Section $3$, we define the
 $(G,m)$-multiparking complement function and study its
properties. In Section $4$, by a combinatorial method, we obtain a
recursion of the generating function of the sum
$\sum\limits_{i=1}^na_i$ of $G$-parking functions
$(a_1,\ldots,a_n)$.
\section{$(G,m)$-multiparking functions}
In this section, first we give the definition of
$(G,m)$-multiparking function. Given a $m\in \{1,2,\ldots,n\}$, for
any $I\subseteq [n]$, let $\alpha(I,m)=\min\{i\in I\mid i\geq m\}$.
\begin{definition}\label{definition} Let $1\leq m\leq n$ and $G$ be a connected graph with vertex set
$V(G)=\{1,2,\cdots,n\}$. A $(G,m)$-multiparking function is a
function $f:V(G)\rightarrow \mathbb{N}\cup\{-1\}$, such that for
every $I\subseteq V(G)$ either (A) $f(\alpha(I,m))=-1$, or (B) there
exists a vertex $i\in I$ such that $0\leq f(i)<outdeg_I(i)$.
\end{definition}
The vertices which satisfy $f(i)=-1$ in (A) will be called roots of
$f$. Furthermore, we say that the vertex $v$ is called a {\it
absolute root} if $f(v)=-1$ in all $(G,m)$-multiparking functions
$f$; and the vertex $v$ is called a {\it relative root} if there are
$(G,m)$-multiparking functions $f$ and $f'$ such that $f(v)=-1$ and
$f'(v)\geq 0$ respectively. By Definition \ref{definition}, it is
easy to see that a $(G,m)$-multiparking function is a
$G$-multiparking functions and a $G$-parking functions when $m=1$
and $m=n$ respectively.

Let $G$ be a connected graph with vertex set $V(G)=[n]$, where $[n]$
denote the set $\{1,2,\ldots,n\}$. There are loops and multiple
edges in $G$. Given $m\in [n]$, let $m$ be the absolute root of $G$.
For any $i,j\in[n]$, let $\mu_G(i,j)$ be the number of edges between
the vertices $i$ and $j$ in $G$. For establishing the bijections,
all edges of $G$ are colored. The colors of edges connecting the
vertices $i$ and $j$ are $0,1,\cdots,\mu_G(i,j)-1$ respectively for
any $i,j\in V(G)$. We use $\{i,j\}_k$ to denote the edge $e\in E(G)$
connecting two vertices $i$ and $j$ with color $k$. A color
$m$-subforest $F$ of $G$ is a color subgraph of $G$ without cycles
such that there is a vertex $j\in [m,n]$ in every component of $F$.
Let $\mathcal{MP}_{G,m}$ and $\mathcal{F}_{G,m}$ be the sets of the
$(G,m)$-multiparking functions and the spanning color $m$-forests of
$G$ respectively. For any $F\in\mathcal{F}_{G,m}$ and $e\in F$, let
$c_F(e)$ denote the color of edge $e$ in $F$. By modifying the
algorithms A and B in \cite{kostic}, we construct the bijection
$\Phi$ between $\mathcal{MP}_{G,m}$ and $\mathcal{F}_{G,m}$. Since
the proof of the bijection is similar to the proof in \cite{kostic},
we only give the sketch of the proof in this paper. The following
algorithm gives a mapping $\Phi$ from $\mathcal{MP}_{G,m}$ to
$\mathcal{F}_{G,m}$.

\noindent{\bf Algorithm A. (Kostic, Yan \cite{kostic})}

{\bf Step 1:} Let $val_0=f$, $P_0=\emptyset$, $F_0=Q_0=\{m\}$.

{\bf Step 2:} At time $i\geq 1$, let $v=\min \{\tau(w)\mid w\in
Q_{i-1}\}$, where $\tau$ is a vertex ranking in $S_n$.

{\bf Step 3:} Let $N=\{w\notin P_{i-1}\mid 0\leq val_i(w)\leq
\mu(w,v)-1\text{ and }\{w,v\}_{val_i(w)}\in E(G)\}$ and
$\hat{N}=\{w\notin P_{i-1}\mid val_i(w)\geq \mu(w,v)\text{ and
}\{w,v\}_{val_i(w)}\in E(G)\}$. Set
$val_i(w)=val_{i-1}(w)-\mu_G(w,v)$ for all $w\in \hat{N}$. For any
other vertex $u$, set $val_i(w)=val_{i-1}(w)$. Update $P_i$, $Q_i$
and $F_i$ by letting $P_i=P_{i-1}\cup\{v\}$, $Q_i=Q_{i-1}\cup
N\setminus \{v\}$ if $Q_{i-1}\cup N\setminus \{v\}\neq \emptyset$,
otherwise $Q_{i}=\{u\}$ where $u=\alpha([n]\setminus P_i,m)$. Let
$F_i$ be a graph on $P_{i}\cup Q_i$ whose edges are obtained from
those of $F_{i-1}$ by joining edges $\{w,v\}_{val_{i-1}(w)}$ for
each $w\in N$.

Define $\Phi=\Phi_{G,m,\tau}:\mathcal{MP}_{G,m}\rightarrow
\mathcal{F}_{G,m}$ by letting $\Phi(f)=F_n$. First, iterating the
Steps $2$-$3$ until $n$, we must have $P_n=[n]$ and $Q_n=\emptyset$.
Otherwise, let $\tilde{P}=[n]\setminus P_n$. Then $f(i)\geq
outdeg_{\tilde{P}}(i)$ for all $i\in \tilde{P}$, a contradiction.
Also, it is easy to see that each $F_i$ is a forest since every edge
$\{v,w\}_k$ in $F_i\setminus F_{i-1}$ has one endpoint in
$V(F_i)\setminus V(F_{i-1})$. In the above algorithm, let $f$ be a
$(G,m)$-multiparking function. Then each tree component $T$ of
$\Phi(f)$ has exactly one vertex $v$ with $f(v)=-1$. In particular,
$v=\alpha(V(T),m)$. Thus, we have
$\Phi(\mathcal{MP}_{G,m})\subseteq\mathcal{F}_{G,m}$.

Let $F\in \mathcal{F}_{G,m}$. Suppose $T$ is a component of $F$. Let
$\alpha(V(T),m)$ be the root of $T$.  For any non-root $v\in [n]$,
there is an unique root $r_i$ which is connected with $v$. Define
the height of $v$ to be the number of edges in the path connecting
$v$ with root $r_i$. If the height of a vertex $w$ is less than the
height of $v$ and $\{v,w\}_k$ is an edge of $F$, then $w$ is the
predecessor of $v$, $v$ is a child of $w$, and write $w={ pre}_F(v)$
and $v\in { child}_F(w)$. The following algorithm will give the
inverse map of $\Phi$.

Let $G$ be a connected color graph with a spanning color
$m$-subforest $F$. A leaf of $F$ is a vertex $v\in V(F)$ with degree
$1$ in $F$. Denote the set of leaves of $F$ by $Leaf(F)$. Let
$T_1,T_2,\ldots T_k$ be the trees of $F$ with respective roots
$m=r_1<r_2<\ldots r_k$, where $r_i=\alpha(V(T_i),m)$.\\
\noindent{\bf Algorithm B (Kostic, Yan \cite{kostic}).}

{\bf Step 1. } Let $\tau$ be a vertex ranking in $S_n$. Assume
$v_1,v_2,\ldots, v_i$ are determined, where $v_1=m$. Let
$V_i=\{v_1,v_2,\ldots, v_i\}$ and $W_i=\{v\notin V_i\mid \{v,w\}\in
F\text{ for some }w\in V_i\}$.

{\bf Case (1)} If $W_i=\emptyset$, let $v_{i+1}$ be the minimum
vertex which is larger than or equal to $m$ in the set
$V(G)\setminus V_i$;

{\bf Case (2)} Otherwise, let $F'$ be the forest obtained by
restricting $F$ to $V_i\cup W$. Let $v_{i+1}=min\{\tau(w)\mid w\in
Leaf(F')\}$.

{\bf Step 2. } Use $\sigma(F)$ to denote the number of the connected
components in $F$. Set $f(r_1)=f(r_2)=\cdots =f(r_{\sigma(F)})=-1$.
For any other vertex $v$, let $f(v)$ be equal to the sum of
 the color of the edge connecting the vertices $v$ with $pre_F(v)$ and the cardinality of the set
 $N(v)$, where
$N(v)=N_{G,F,\tau}(v)=\{v_j\mid \{v,v_j\}_k\in E(G)\text{ and }
\pi^{-1}(v_j)<\pi^{-1}({\rm pre}_F(v))\}$.

Define $\Psi=\Psi_{G,m,\tau}:\mathcal{F}_{G,m}\rightarrow
\mathcal{MP}_{G,m}$ by letting $\Psi_{G,m,\tau}(f)=f_F$. For any
$I\subseteq [n]$, Let $\pi$ be the permutation defined in Step 1 of
the Algorithm B. It is easy to see that there exists an unique
integer $k$ such that $I\subseteq \{\pi(k),\ldots,\pi(n)\}=I'$ and
$\pi(k)\in I$. If the vertex $\pi(k)$ isn't the root of a component
in $F$, then we have $0\leq f(\pi(k))<outdeg_{I'}(\pi(k))$. So,
$f(\pi(k))<outdeg_{I}(\pi(k))$ since $outdeg_{I}(\pi(k))\leq
outdeg_{I'}(\pi(k))$. If the vertex $\pi(k)$ is the root of a
component in $F$, then $f(\pi(k))=-1$ and $\pi(k)=\alpha(I',m)$.
Clearly, $\pi(k)=\alpha(I,m)$ since $I\subseteq I'$. Hence, by the
Definition \ref{definition}, we have $f=f_F$ is a
$(G,m)$-multiparking function and
$\Psi(\mathcal{F}_{G,m})\subseteq\mathcal{MP}_{G,m}$.

Finally, note that the order $\pi=v_1v_2\ldots v_n$ in the algorithm
B is exactly the order in which vertices of $G$ will be placed into
the set $P_i$ when running algorithm A on $f$. Thus, we have
$\Psi(\Phi(f))=f$. Hence, $\Phi=\Phi_{G,m,\tau}$ and
$\Psi=\Psi_{G,m,\tau}$ are inverses of each other. We state these
results as the following theorem.
\begin{thm}\label{theoremG,mmultiparkingfunction}(Kostic, Yan \cite{kostic}) The mapping $\Phi$ is a bijection from $\mathcal{MP}_{G,m}$ to
 $\mathcal{F}_{G,m}$.
\end{thm}

By Algorithm B, we immediately obtain the following corollary.
\begin{cor}\label{corollaryfF=c(e)+N-sigma}(Kostic, Yan \cite{kostic}) For any $F\in\mathcal{F}_{G,m}$, let $f_F$ be the
$(G,m)$-multiparking function  corresponding with $F$, i.e.,
$f_F=\Phi^{-1}(F)$. Let $\sigma(F)$ be the number of the connected
components in $F$ and $R(F)$ the set of all the roots in $F$. Then
$\sum\limits_{v\in[n]}f_F(v)=\sum\limits_{e\in
E(F)}c(e)+\sum\limits_{v\in V(G)\setminus R(F)}|N(v)|-\sigma(F)$.
\end{cor}
\begin{proof} For any $F\in \mathcal{F}_{G,m},
$ let $f_F=\Phi^{-1}(F)$. Recall that $\sigma(F)$ is the number of
the connected components in $F$. For any $v\in [n]$, if $v$ is  the
minimal vertex which is no less than $m$ in the tree of $F$, then
$f_F(v)=-1$; otherwise, $f_F(v)$ is the sum of
 the color of the edge $e$ connecting the vertices $v$ to $pre_F(v)$ and the cardinality of the set $N(v)$. Hence,
\begin{eqnarray*}\sum\limits_{v\in[n]}f_F(v)&=&\sum\limits_{f_F(v)\neq -1}f_F(v)+\sum\limits_{f_F(v)=-1}f_F(v)\\
&=&\sum\limits_{f_F(v)\neq
-1}f_F(v)-\sigma(F)\\
&=&\sum\limits_{e\in E(F)}c_F(e)+\sum\limits_{v\in V(G)\setminus
R(F)}|N(v)|-\sigma(F)\\\end{eqnarray*}
\end{proof}
\section{$(G,m)$-multiparking
complement  functions} In this section, we will define  the
$(G,m)$-multiparking complement function, give the reciprocity
theorem for $(G,m)$-multiparking function. First, let $deg_G(i)$ be
the number of edges which are incident with the vertex $i$ in the
graph $G$. For any subset $I\subseteq V(G)$ and $i\in I$, define
$indeg_I(i)$ as the number of edges from $i$ to vertices inside $I$.
We give the definition of $(G,m)$-multiparking complement  function
as follows.
\begin{definition} Let $m\in [n]$ and $G$ be a connected graph with vertex set
$V(G)=[n]$. A
 $(G,m)$-multiparking complement function is a function
$h:V(G)\rightarrow \mathbb{N}$, such that for every $I\subseteq
V(G)$ either ($\overline{A}$) $h(\alpha(I,m))=deg_G(\alpha(I,m))+1$,
or ($\overline{B}$) there exists a vertex $i\in I$ such that
$indeg_I(i)< h(i)\leq deg_G(i)$.
\end{definition}

Given a function $f:V(G)\rightarrow \mathbb{N}\cup\{-1\}$, define a
function $h_f$ as $h_f(i)={\rm deg}_G(i)-f(i)$ for all $i\in V(G)$.
The following lemma tells us the relation between
$(G,m)$-multiparking functions and $(G,m)$-multiparking complement
functions.
\begin{lemma}  Let $m\in [n]$ and $G$ be a connected graph with vertex set
$V(G)=[n]$.  Then the function $f$ is a $(G,m)$-multiparking
function if and only if the function $h_f$ is a $(G,m)$-multiparking
complement function.
\end{lemma}
\begin{proof} Let $f\in\mathcal{MP}_{G,m}$. For any $I\subseteq V(G)$, if
$I$ contains no well-behaved vertices, then $f(\alpha(I,m))=-1$,
thus, $h_f(\alpha(I,m))={\rm deg}_G(\alpha(I,m))+1$; otherwise,
there exists a vertex $i\in I$ such that $0\leq f(i)< {\rm
outdeg}_I(i)$, then ${\rm indeg}_I(i)<h_f(i)\leq {\rm deg}_{G}(i)$.
Hence, $h_f$ is a $(G,m)$-multiparking complement function.
\end{proof}

Let $\overline{\mathcal{MP}}_{G,m}$ be a set of all the
$(G,m)$-multiparking complement functions. The generating function
$P_{G,m}(q)$ of $(G,m)$-multiparking functions is defined as
\begin{eqnarray*}P_{G,m}(q)=\sum\limits_{f\in \mathcal{MP}_{G,m}}q^{\sum\limits_{i\in V(G)}f(i)}\end{eqnarray*}
Define the generating function of $(G,m)$-multiparking complement
functions to be the
polynomial\begin{eqnarray*}\overline{P}_{G,m}(q)&=&\sum\limits_{h\in
\overline{\mathcal{MP}}_{G,m}}q^{\sum\limits_{i\in
V(G)}h(i)-|E(G)|}\\
&=&\sum\limits_{f\in {\mathcal{MP}}_{G,m}}q^{\sum\limits_{i\in
V(G)}[deg_G(i)-f(i)]-|E(G)|}\\
&=&q^{|E(G)|}P_{G,m}\left(\frac{1}{q}\right)\end{eqnarray*} Clearly,
the sums $\sum\limits_{i\in V(G)}[deg_G(i)-f(i)]-|E(G)|$ and
$\sum\limits_{i\in V(G)}f(i)$ are connected with the reciprocity law
for $(G,m)$-mulitiparking function. Now, we are in a position to
give the reciprocity theorem for $(G,m)$-multiparking function.

A color forest $F$ on $[n]$ may appear as a subgraph of different
graphs, and a vertex function $f$ may be a $(G,m)$-multiparking
function for different graphs. Let $f=\Phi(F)$. Recall that $G-e$ be
a graph obtained by deleting the edge $e=\{i,j\}_k$ from $G$. We say
an edge $e=\{i,j\}_k$ of $G-F$ is $F$-{\it redundant} if
$\Phi_{G-e,m}(F)=f$. There are the closed relations between
$(G,m)$-multiparking complement functions and $F$-redundant edges of
$G$.

Let $\pi$ be the order defined in Step 1 of Algorithm $B$. Note that
$\pi$ only depends on $F$ and $\tau$, not the underlying graph $G$.
We have the following lemma.
\begin{lemma}\label{redundant}(Kostic, Yan \cite{kostic}) Fix a vertex ranking $\tau$. An edge $e=\{v,w\}_k$ of $G$ is $F$-redundant if and
only if $e$ is one of the following types:

1. Both $v$ and $w$ are roots of $F$.

2. $v$ is a root and $w$ are non-roots of $F$, and
$\pi^{-1}(w)<\pi^{-1}(v)$.

3. $v$ and $w$ are non-roots and
$\pi^{-1}(pre_F(v))<\pi^{-1}(w)<\pi^{-1}(v)$. In this case $v$ and
$w$ must lie in the same tree of $F$.

4. $e$ is a loop of $G$.

5. $e=\{v,pre_F(v)\}_k$ with $k>k'$, where $\{v,pre_F(v)\}_{k'}\in
F$.
\end{lemma}

The results $1$, $2$ and $3$ in Lemma \ref{redundant} were proved in
\cite{kostic}. It is easy to see the results $4$ and $5$ hold by the
Algorithm $A$.

 By the above lemma, given a graph $G$, let $F$ be a spanning color $m$-subforest of
 $G$. Define a function $g_F:V(G)\rightarrow \mathbb{N}$ by letting
 $g_F(i)$ to be the cardinality of the set $\{e=\{i,j\}_k\mid e \text{ is F-redundant and
 }\pi^{-1}(j)\leq\pi^{-1}(i)\}$. Define the generating function  to be the
polynomial\begin{eqnarray*}I_{G,m}(q)&=&\sum\limits_{F\in
\mathcal{F}_{G,m}}q^{\sum\limits_{i\in V(G)}g_F(i)}\end{eqnarray*}

\begin{theorem}(the reciprocity
theorem for $(G,m)$-multiparking
function)\begin{eqnarray*}q^{|V(G)|}I_{G,m}(q)=\overline{P}_{G,m}(q)=q^{|E(G)|}P_{G,m}\left(\frac{1}{q}\right).\end{eqnarray*}
\end{theorem}
\begin{proof} We only prove the first identity. Note that $g_F(i)$ is the
cardinality of the set $\{e=\{i,j\}_k\mid e \text{ is F-redundant
and
 }\pi^{-1}(j)\leq\pi^{-1}(i)\}$. By Corollary \ref{corollaryfF=c(e)+N-sigma}, we have \begin{eqnarray*}|E(G)|&=&\sum\limits_{i\in V(G)}g_F(i)+\sum\limits_{f_F(i)\neq -1}|N(i)|+\sum\limits_{e\in
F}c_F(e)+|E(F)|\\
&=&\sum\limits_{i\in V(G)}g_F(i)+\sum\limits_{f_F(i)\neq
-1}|N(i)|+\sum\limits_{e\in
F}c_F(e)+|V(G)|-\sigma(F)\\
&=&\sum\limits_{i\in V(G)}g_F(i)+\sum\limits_{i\in
V(G)}f_F(i)+|V(G)|.\end{eqnarray*} Therefore,
\begin{eqnarray*}\overline{P}_{G,m}(q)&=&\sum\limits_{f\in
{\mathcal{MP}}_{G,m}}q^{|E(G)|-\sum\limits_{i\in
V(G)}f(i)}\\
&=&\sum\limits_{F\in\mathcal{F}_{G,m}}q^{\sum\limits_{i\in V(G)}g_F(i)+V(G)}\\
&=&q^{|V(G)|}I_{G,m}(q).\end{eqnarray*}\end{proof}

\section{The recursion for the generating functions $P_{G,n}(q)$}
In this section, we always let $m=n$. In this situation, a spanning
color $n$-subforest of $G$ is a spanning color  tree of $G$ with
root $n$ since $G$ is connected, and a $(G,n)$-multiparking function
is a $G$-parking function. Let $\mathcal{T}_G=\{T\mid T\text{ is a
spanning color tree of }G\}$ and $\mathcal{P}_{G}=\{f\mid f\text{ is
a }G\text {-parking function}\}$. We write $P_{G,n}(q)$ as
$P_{G}(q)$ for short. There is a bijection $\Phi=\Phi_{G,\tau}$
between the sets $\mathcal{T}_G$ and $\mathcal{P}_G$. For any
$T\in\mathcal{T}_G$, define
$w_G(T)=\sum\limits_{i=1}^{n}\Phi(T)(i)$. Then
$P_{G}(q)=\sum\limits_{T\in \mathcal{T}_G}q^{w_G(T)}$.

First, the Algorithms $A$ and $B$ tell us that $G$-parking functions
are independent on the loops in $G$. So, we obtain the following
result.

\begin{lemma}\label{loop} Suppose that $e$ is a loop of $G$. Then
$P_{G}(q)=P_{G-e}(q)$.
\end{lemma}

Now, suppose that $G$ has a bridge $e=\{i,j\}_0$. After deleting the
edge $e$, we let $G_1$ be the subgraph of $G$ such that
$\{n,i\}\subseteq V(G_1)$ and $G_2$ another subgraph of $G$ obtained
by letting the label of the vertex $j$ become $n+1$.

\begin{lemma}\label{bridge} Suppose that $e=\{i,j\}_0$ is a bridge of $G$.  Then $P_{G}(q)=qP_{G_1}(q)P_{G_2}(q)$.
\end{lemma}
\begin{proof}
For any $T_1\in\mathcal{T}_{G_1}$ and $T_2\in\mathcal{T}_{G_2}$, the
trees $T_1$ and $T_2$ have the roots $n$ and $n+1$ respectively. Let
$T$ be a tree obtained by setting the label of the vertex $n+1$ as
$j$ and adding $j$ to be the child of the vertex $i$. Then $T$ is a
spanning color tree of $G$. Conversely, for any
$T\in\mathcal{T}_{G}$, after deleting the edge $e$, we let $T_1$ be
the subgraph of $T$ such that $\{n,i\}\subseteq V(T_1)$ and $T_2$
another subgraph of $T$ obtained by letting the label of the vertex
$j$ become $n+1$. Then $T_i\in \mathcal{T}_{G_i}$ for $i=1,2$.  Take
the vertex ranking $\tau$ such that $\tau(i)=n-1$ and $\tau(j)=n-2$.
Then $w_{G}(T)=w_{G_1}(T_1)+ w_{G_2}(T_2)+1$ since $e$ is a bridge
of $G$. This implies that
\begin{eqnarray*}\sum\limits_{T\in\mathcal{T}_{G}}q^{w_{G}(T)}=q\sum\limits_{T\in\mathcal{T}_{G_1}}q^{w_{G_1}(T)}
\sum\limits_{T\in\mathcal{T}_{G_2}}q^{w_{G_2}(T)}.\end{eqnarray*}
Hence, we have $P_{G}(q)=qP_{G_1}(q)P_{G_2}(q)$.
\end{proof}

Now, we consider the case in which $e=\{i,j\}_k$ is neither a loop
nor a bridge with $i>j$. Define a graph $G{\setminus e}$ as follows.
The graph $G{\setminus e}$ is obtained from $G$ contracting the the
vertices $i$ and $j$; that is, to get $G{\setminus e}$ we identify
two vertices $i$ and $j$ as a new vertex $i$. First, we discuss the
cases in which $e$ is incident with the root $n$.

\begin{lemma}\label{lemmae=n1}Suppose that $e=\{n,i\}_k\in E(G)$ and $e$ is neither a loop nor a bridge in
$G$. Then $P_{G}(q)=qP_{G-e}(q)+P_{G\setminus e}(q)$.
\end{lemma}
\begin{proof} First, we take the vertex ranking $\tau$ such that $\tau(i)=1$.
Let $\mathcal{T}_{G,0}=\{T\in\mathcal{T}_G\mid e\notin T\}
\cup\{T\in\mathcal{T}_G\mid e\in T\text{ and }c_T(e)\geq 1\}$ and
$\mathcal{T}_{G,1}=\{T\in\mathcal{T}_G\mid e\in T\text{ and
}c_T(e)=0\}$. Given a spanning color tree $T$ of $G$, let
$\pi=\pi_{T,\tau}$ be a permutation $\pi=(\pi(1),\pi(2),\ldots,
\pi(n))=(v_1v_2\ldots v_n)$ on the vertices of $G$ by Step $1$ of
the Algorithm B. For any $T\in\mathcal{T}_{G,0}$, if $e\notin T$,
suppose $j=pre_{T}(i)$, then $j\neq n$ and
$\pi^{-1}(j)+1=\pi^{-1}(i)$ since $\tau(i)=1$. Thus
$w_{G-e}(T)=w_{G}(T)-1$ since $T$
 is  a spanning color tree of $G-e$ as well. If $e\in T$ and $c_T(e)\geq
1$, then let $T'$ be the tree obtained by setting the color of the
edge $e$ as $c_T(e)-1$. The tree $T'$ is a spanning color tree of
$G-e$ and $w_{G-e}(T)=w_{G}(T)-1$.

For any $T\in\mathcal{T}_{G,1}$, it is easy to see that
$pre_{T}(i)=n$ and $\pi^{-1}(i)=2$. Let $A$ be the set of the
vertices $j\in {\rm child_T}(i)$ such that $ n\in N(j)$. Since
$c_T(e)=0$, let $T'$ be the tree obtained by deleting the vertex
$i$, attaching the vertex $j$ to be the child of $n$ and setting the
color of the edge $\{n,j\}$ as $c_T(\{i,j\})+\mu_G(n,j)$ for all
$j\in A$. The tree $T'$ can be viewed as a spanning color tree of
$G\setminus e$. Clearly, $w_{G\setminus e}(T')=w_G(T)$. Conversely,
Given a tree $T'\in \mathcal{T}_{G\setminus e}$,  let $A$ be the set
of the vertices which  is adjacent to $n$ in $T$ such that either
$c_T(\{n,j\})\geq \mu_G(n,j)$ or $j$ isn't adjacent to $n$ in $G$.
Add a new vertex $i$, attach $i$ to be the child of $n$ and let the
color of the edge $\{n,i\}$ be $0$. Then delete the edges $\{n,j\}$,
attach $j$ to be the child of $i$ and setting the color of $\{n,j\}$
as $c_T(\{n,j\})-\mu_G(n,j)$ for all $j\in A$. The obtained tree can
be viewed as a spanning color tree of $G$.
 Hence, $P_{G}(q)=qP_{G-e}(q)+P_{G\setminus e}(q)$.
\end{proof}

Next, we discuss the case in which $e$ isn't incident with the root
$n$. Suppose $e=\{i,j\}_k\in E(G)$ with $i>j$. Take the vertex
ranking $\tau$ such that $\tau(i)=1$ and $\tau(j)=2$. Define the
following six sets:

(1) $\hat{T}_{G,0}^1=\{T\in\mathcal{T}_{G}\mid e\notin T\text{ and
}pre_{T}(j)\neq pre_{T}(i)\}$

(2) $\hat{T}_{G,0}^2=\{T\in\mathcal{T}_{G}\mid e\notin T,
pre_{T}(j)= pre_{T}(i)\text{ and }c_T(\{j,pre_{T}(i)\})\neq 0\}$

(3) $\hat{T}_{G,0}^3=\{T\in\mathcal{T}_{G}\mid e\notin T,
pre_{T}(j)= pre_{T}(i)\text{ and }c_T(\{j,pre_{T}(i)\})= 0\}$

(4) $\hat{T}_{G,1}^1=\{T\in\mathcal{T}_{G}\mid e\in T\text{ and
}c_T(e)\geq 1\}$

(5) $\hat{T}_{G,1}^2=\{T\in\mathcal{T}_{G}\mid e\in T,
\pi^{-1}_{T}(i)<\pi^{-1}_{T}(j), pre_{T}(i)\in N_{G,T}(j)\text{ and
}c_T(e)=0\}$, where $N_{G,T}(v)=\{v_j\mid \{v,v_j\}_k\in E(G)\text{
and } \pi^{-1}(v_j)<\pi^{-1}({\rm pre}_F(v))\}$.

(6) $\hat{T}_{G,1}^3=\{T\in\mathcal{T}_{G}\mid e\in
T\}\setminus\left(\hat{T}_{G,1}^1\cup\hat{T}_{G,2}^2\right)$.

\begin{lemma}\label{lemmarecurrsionni}Suppose that $e=\{i,j\}_k\in E(G)$ and $e$ is neither a loop nor a  bridge in $G$.
There is a bijection $\psi$ from
$\hat{T}_{G,0}^1\cup\hat{T}_{G,0}^2\cup\hat{T}_{G,1}^1\cup\hat{T}_{G,1}^2$
to $T_{G-e}$. Moreover, $w_{G-e}(\psi(T))=w_{G}(T)-1$ for any $T\in
\hat{T}_{G,0}^1\cup\hat{T}_{G,0}^2\cup\hat{T}_{G,1}^1\cup\hat{T}_{G,1}^2$.
\end{lemma}
\begin{proof} For any $T\in\hat{T}_{G,0}^1$, let $\psi(T)=T$. $\psi(T)$ can
be viewed as a spanning tree of $G-e$. Moreover,
$w_{G-e}(T)=w_{G}(T)-1$.

For any $T\in\hat{T}_{G,0}^2$, suppose $v=pre_T(i)$. Since
$c_T(\{j,v\})\neq 0$, let $T'$ be the tree obtained by setting the
color of $\{j,v\}$ as $c_T(\{j,v\})-1$. The tree $T'$ can be viewed
as a spanning tree of $G-e$ and satisfies that
$pre_{T'}(i)=pre_{T'}(j)$ and $c_{T'}(\{j,v\})\leq
\mu_{G-e}(j,v)-2=\mu_{G}(j,v)-2$. Moreover, $w_{G-e}(T')=w_G(T)-1$.

For any $T\in\hat{T}_{G,1}^1$, let $T'$ be the tree obtained by
setting the color of $\{i,j\}$ as $c_T(i,j)-1$. The tree $T'$ can be
viewed as a spanning tree of $G-e$ and $c_{T'}(\{i,j\})\geq 0$.
Moreover, $w_{G-e}(T')=w_G(T)-1$.

For any $T\in\hat{T}_{G,1}^2$, suppose that $v=pre_{T}(i)$. Let $T'$
be the tree obtained by deleting the edge $\{i,j\}_0$, attaching $j$
to be the child of $v$ and setting the color of  the edge $\{v,j\}$
as $\mu_G(v,j)-1$. Since $v\in N(j)$, the tree $T'$ can be viewed as
a spanning tree of $G-e$,
$c_{T'}(\{j,v\})=\mu_{G-e}(j,v)-1=\mu_{G}(j,v)-1$ and
$pre_{T'}(i)=pre_{T'}(j)$. Moreover, $w_{G-e}(T')=w_G(T)-1$. This
complete the proof.\end{proof}

\begin{lemma}\label{lemmarecurrsionij}Suppose that $e=\{i,j\}_k\in E(G)$ and $e$ is neither a loop nor a  bridge in $G$.
There is a bijection $\psi'$ from
$\hat{T}_{G,0}^3\cup\hat{T}_{G,1}^3$ to $T_{G\setminus e}$.
Moreover, $w_{G\setminus e}(\psi(T))=w_{G}(T)$ for any
$\hat{T}_{G,0}^3\cup\hat{T}_{G,1}^3$.
\end{lemma}
\begin{proof} For any $T\in\hat{T}_{G,0}^3$, suppose $v=pre_T(j)$. Let $T'$
be the tree obtained by deleting the vertex $j$, attaching the
vertex $w$ to be the child of $i$ and setting the color of $\{i,w\}$
as $c_T(\{j,w\})+\mu_G(i,w)$ for any $w\in {\rm child}_T(j)$. The
tree $T'$ can be viewed as a spanning tree of $G\setminus e$ and
$0\leq c_{T'}(\{i,v\})\leq \mu_G(\{i,v\})-1$. Moreover,
$w_{G\setminus e}(T')=w_{G}(T)$.

For any $T\in\hat{T}_{G,1}^3$, we always have $c_T(\{i,j\})=0$. we
discuss the following three cases:

{\it Case 1.} $\pi_T^{-1}(i)<\pi_T^{-1}(j)$ and $pre_{T}(i)\notin
N(j)$

Let $T'$ be the tree obtained by deleting the vertex $j$, attaching
$w$ to be the child of $i$ and setting the color of the edge
$\{w,i\}$ as $c_T(\{j,w\})+\mu_G(i,w)$ for all $w\in child_T(j)$.
The tree $T'$ can be viewed as a spanning tree of $G\setminus e$.
Moreover, $w_{G\setminus e}(T')=w_G(T)$.

{\it Case 2.} $\pi_T^{-1}(i)>\pi_T^{-1}(j)$ and $pre_{T}(j)\notin
N(i)$

Suppose $v=pre_{T}(j)$. Let $T'$ be the tree obtained by deleting
the vertex $j$, attaching $i$ to be the child of $v$, attaching $w$
to be the children of $i$ for all $w\in{\rm child}_T(j)$  and
setting the color of the edge $\{u,i\}$ as $c_T(\{i,u\})+\mu_G(j,u)$
for all $u\in {\rm child}_T(i)$. The tree $T'$ can be viewed as a
spanning tree of $G\setminus e$. Moreover, $w_{G\setminus
e}(T')=w_G(T)$.

 {\it Case 3.}
$\pi_T^{-1}(i)>\pi_T^{-1}(j)$ and $pre_{T}(j)\in N_{G,T}(i)$

Suppose that $v=pre_{T}(j)$. Let $T'$ be the tree obtained by
deleting the vertex $j$, attaching $i$ to be the child of $v$,
attaching $w$ to be the children of $i$ for all $w\in{\rm
child}_T(j)$ and setting the color of the edge $\{w,i\}$ as
$c_T(\{i,w\})+\mu_G(j,w)$ for all $w\in child_T(i)$ and the color of
the edge $\{v,i\}$ as $c_T(\{j,v\})+\mu_G(i,v)$. The tree $T'$ can
be viewed as a spanning tree of $G\setminus e$ and
$c_{T'}(\{i,v\})\geq \mu_{G}(i,v)$. Moreover, $w_{G-e}(T')=w_G(T)$.

Conversely, for any $T\in \mathcal{T}_{G\setminus e}$, suppose
$v=pre_{T}(i)$. we consider the following three cases.

{\it Case 1.} $\mu_G(i,v)\geq 1$, $\mu_G(j,v)\geq 1$ and $0\leq
c_{T}(\{i,v\})\leq\mu_G(i,v)-1$.

Add a new vertex $j$ and attach $j$ to be the child of $v$. Let the
color of $\{v,j\}$ be $0$. For any $w\in {\rm child}_T(i)$, if
$c_T(\{i,w\})\geq \mu_{G}(i,w)$ or $\mu_G(i,w)=0$, then delete the
edge $\{i,w\}$ and attach $w$ to be the child of $j$, let the color
of $\{j,w\}$ be $c_T(\{i,w\})-\mu_G(i,w)$. The obtained tree can be
view as a spanning tree of $G$.

{\it Case 2.} $\mu_G(i,v)\geq 1$, $\mu_G(j,v)\geq 1$ and
 $c_{T}(\{i,v\})\geq \mu_G(i,v) $

Add a new vertex $j$ and attach $j$ to be the child of $v$. Delete
the edge $\{v,i\}$ and attach $i$ to be the child of $j$. Let the
color of $\{v,j\}$ be $c_T(\{i,v\})-\mu_G(i,v)$ and the color of
$\{i,j\}$ 0. For any $w\in {\rm child}_T(i)$, if $c_T(\{i,w\})\geq
\mu_{G}(j,w)$ or $\mu_G\{j,w\}=0$, then let the color of the edges
$\{i,w\}$ be $c_T(\{i,w\})-\mu_G(j,w)$; otherwise delete the edge
$\{i,w\}$ and attach $w$ to be the child of $j$. The obtained tree
can be view as a spanning tree of $G$.

{\it Case 3.} There is exact one vertex $w$ such that
$\mu_G(w,v)\geq 1$ for $w\in \{i,j\}$

We first consider the case $\mu_G(i,v)\geq 1$, $\mu_G(j,v)=0$. Add a
new vertex $j$ and attach $j$ to be the child of $i$. Let the color
of $\{i,j\}$ be $0$. For any $w\in {\rm child}_T(i)$, if
$c_T(\{i,w\})\geq \mu_{G}(i,w)$  or $\mu_G(i,w)=0$, then delete the
edge $\{i,w\}$ and attach $w$ to be the son of $j$.  Let the color
of $\{j,w\}$ be $c_T(\{i,w\})-\mu_G(i,w)$. The obtained tree can be
view as a spanning tree of $G$. Similarly, we may consider the case
$\mu_G(i,v)=0$, $\mu_G(j,v)\geq 1$. The proof is completed.
\end{proof}

\begin{theorem}\label{noloopnobridge} Suppose that $e\in E(G)$ and $e$ is neither a loop nor a  bridge in $G$. Then
$$P_{G}(q)= qP_{G-e}(q)+P_{G\setminus e}(q).$$
\end{theorem}
\begin{proof} Combining Lemmas \ref{lemmarecurrsionni} and
\ref{lemmarecurrsionij}, we obtain the results as desired.
\end{proof}


\end{document}